\documentclass[11pt,reqno]{amsart}
\usepackage{amsmath,amssymb,graphicx}
\newtheorem{thm}{Theorem}[section]
\newtheorem{cor}[thm]{Corollary}
\newtheorem{lem}[thm]{Lemma}
\newtheorem{prop}[thm]{Proposition}
\theoremstyle{definition}
\newtheorem{defn}[thm]{Definition}
\theoremstyle{remark}
\newtheorem{rem}[thm]{Remark}
\theoremstyle{remark}
\newtheorem{ex}[thm]{Example}
\numberwithin{equation}{section}

\begin{document}

\title[Multiplicative decompositions]{Multiplicative decompositions and frequency of vanishing of nonnegative submartingales}%
\author{Ashkan Nikeghbali}
\address{American Institute of Mathematics
360 Portage Ave Palo Alto, CA 94306-2244 \\
and University of Rochester} \email{ashkan@aimath.org}
 \subjclass[2000]{05C38, 15A15;
05A15, 15A18} \keywords{Random Times, submartingales, General theory
of stochastic processes.}
\date{\today}
\begin{abstract}
In this paper, we establish a multiplicative decomposition formula
for nonnegative local martingales and use it to characterize the set
of continuous local submartingales $Y$ of the form $Y=N+A$, where
the measure $dA$ is carried by the set of zeros of $Y$. In
particular, we shall see that in the set of all local submartingales
with the same martingale part in the multiplicative decomposition,
these submartingales are the smallest ones. We also study some
integrability questions in the multiplicative decomposition and
interpret the notion of saturated sets in the light of our results.
\end{abstract}
\maketitle

\section{Introduction}
In \cite{AshkanremI,AshkanremII} we have introduced the following
family of nonnegative local submartingales which enjoy many
remarkable properties, and whose definition goes back to Yor
\cite{yorinegal}:
\begin{defn}\label{martreflechies}
Let $\left(X_{t}\right)$ be a nonnegative local submartingale, which
decomposes as:
\begin{equation}\label{dcodexdbmey}
X_{t}=N_{t}+A_{t}.
\end{equation}
We say that $\left(X_{t}\right)$ is of class $(\Sigma)$ if:
\begin{enumerate}
\item $\left(N_{t}\right)$ is a continuous local martingale, with $N_{0}=0$;
\item $\left(A_{t}\right)$ is a continuous increasing process, with $A_{0}=0$;
\item the measure $\left(dA_{t}\right)$ is carried by the set
$\left\{t:\;X_{t}=0\right\}$.
\end{enumerate}If additionally, $\left(X_{t}\right)$ is of class
$(D)$, we shall say that $\left(X_{t}\right)$ is of class $(\Sigma
D)$.
\end{defn}
A natural question to ask is: what is the place of the local
submartingales of class $(\Sigma)$ among other nonnegative local
submartingales? In particular, is there a systematic way to
construct them? This paper has two aims:
\begin{itemize}
\item to understand better the nature of the
local submartingales of class $(\Sigma)$;
\item to give a characterization which will provide us with further
examples of such processes.
\end{itemize}
Our approach is based on the multiplicative decompositions of
nonnegative local submartingales: indeed, as will be shown, the
stochastic processes of class $(\Sigma)$ reveal some of their very
nice properties when they are factorized (in particular, in a sense
that will be made precise, these processes have more zeros than
other nonnegative local submartingales).

The theory of multiplicative decompositions of semimartingales has
not received much attention in the literature and is much less
powerful than the additive or Doob-Meyer decompositions of
semimartingales. But in the special case of interest to us, it will
be the relevant decomposition to consider.

One of the very first results on multiplicative decompositions, and
which is not so well known, is due to It\^{o} and Watanabe
(\cite{itowatanbe}) who established a multiplicative decomposition
for nonnegative supermartingales in the framework of Markov
Processes: any nonnegative supermartingale $\left(Z_{t}\right)$ is
factorized as:
$$Z_{t}=Z_{t}^{(0)}Z_{t}^{(1)},$$with a positive local martingale
$Z_{t}^{(0)}$ and a decreasing process $Z_{t}^{(1)}$
($Z_{0}^{(1)}=1$).
 Later, Meyer and Yoeurp
(\cite{meyermult},\cite{meyeryoeurp}) and Az\'{e}ma
(\cite{azemamult}) have studied multiplicative decompositions of
nonnegative submartingales and bounded nonnegative supermartingales
to find conditions under which there exists a (unique) raw
increasing process $K$ which generates them. Their approach differed
from that of It\^{o} and Watanabe and was based on the then
flourishing general theory of stochastic processes. The results
obtained in \cite{meyermult,meyeryoeurp,azemamult} are very refined
and too technical for our purpose so we will not discuss them here.
In fact, they do not even apply very well to our setting. Let us
eventually mention that it is proved in \cite{ashyordoob} that every
supermartingale $\mathbb{P}\left(L>t|\mathcal{F}_{t}\right)$
associated with an honest time $L$ can be represented as
$\frac{N_{t}}{\overline{N}_{t}}$, where
$\overline{N}_{t}=\sup_{u\leq t}N_{u}$, and where $N$ is a
nonnegative local martingale starting from one and vanishing at
infinity. This in turn entailed some nice connections with the
theory of enlargements of filtrations. Let us eventually mention the
book of Jacod \cite{jacod} which contains a paragraph about
multiplicative decompositions (but the result of It\^{o} and
Watanabe is not mentioned there).

The main results of this paper are the following:
\begin{enumerate}
\item Every continuous nonnegative local submartingale $\left(Y_{t}\right)$, with $Y_{0}=0$
decomposes uniquely as:
\begin{equation}\label{decomultgenI}
    Y_{t}=M_{t}C_{t}-1,
\end{equation}where $\left(M_{t}\right)$ is a continuous nonnegative local martingale with
$M_{0}=1$ and $\left(C_{t}\right)$ is a continuous increasing
process, with $C_{0}=1$. Moreover, $\left(Y_{t}\right)$ is of class
$(D)$ if and only if $\left(M_{t}\right)$ is a uniformly integrable
martingale and
$$\mathbb{E}\left[M_{\infty}C_{\infty}\right]<\infty.$$
\item A stochastic process $\left(X_{t}\right)$ is of class
$(\Sigma)$ if and only if \begin{equation}\label{decomultgenII}
    X_{t}=\dfrac{M_{t}}{I_{t}}-1,
\end{equation}where $\left(M_{t}\right)$ is a continuous nonnegative local martingale with
$M_{0}=1$ and  $I_{t}=\inf_{u\leq t}M_{u}$. Moreover,
$\left(X_{t}\right)$ is of class $(\Sigma D)$ if and only if
$\left(M_{t}\right)$ is a uniformly integrable martingale and
$$\mathbb{E}\left[\log\dfrac{1}{I_{\infty}}\right]<\infty.$$
\item Let $\left(M_{t}\right)$ be a continuous nonnegative local martingale with
$M_{0}=1$ and let $\mathcal{E}\left(M\right)$ be the set of all
nonnegative local submartingales with the same martingale part $M$
in their multiplicative decomposition. Then, there exists a unique
local submartingale $Y^{\star}$ of class $(\Sigma)$ in
$\mathcal{E}\left(M\right)$ and it is the smallest element in
$\mathcal{E}\left(M\right)$:
$$\forall Y\in\mathcal{E}\left(M\right),\;Y^{\star}\leq Y.$$
\end{enumerate}
This last property is reminiscent of a result of Az\'{e}ma and Yor:
 the set of zeros of submartingales of class $(\Sigma D)$  is a saturated
 set (\cite{azemayorzero}). We shall use our results to give some
 non trivial example of saturated sets $H$ which are not necessarily of
 the form $H=\left\{t:\;M_{t}=0\right\}$ for some martingale $M$.
\section{Multiplicative decompositions}
We first prove the general decomposition formula
(\ref{decomultgenI}), and then specialize to the case when $X$ is of
class $(\Sigma)$. We end the section with a discussion about
saturated sets.
\subsection{General formulae and integrability}
Let $\left( \Omega ,\mathcal{F},\left( \mathcal{F}_{t}\right) _{t\geq 0},%
\mathbb{P}\right) $ a filtered probability space, satisfying the
usual assumptions, and let $Y$ be a continuous nonnegative local
submartingale. It is in general impossible to obtain a
multiplicative decomposition of the form
$Y_{t}=Y_{t}^{(0)}Y_{t}^{(1)},$ with a positive local martingale
$Y_{t}^{(0)}$ and an increasing process $Y_{t}^{(1)}$
($Y_{0}^{(1)}=1$). Indeed, it is well known that once a nonnegative
local martingale is equal to zero, then it remains null, while this
is not true for nonnegative submartingales (consider for example
$|B_{t}|$, the absolute value of a standard Brownian Motion). Hence,
such multiplicative decompositions can hold only for strictly
positive submartingales. But this is not of interest to us since the
local submartingales of the class $(\Sigma)$ have usually many
zeros. That's why we propose the following multiplicative
decomposition.
\begin{prop}
\label{decompositionmultipl}Let $\left(Y_{t}\right)_{t\geq 0}$ be a
continuous nonnegative local submartingale such that $Y_{0}=0$.
Consider its Doob-Meyer decomposition:
\begin{equation}
Y_{t}=m_{t}+\ell_{t}.
\end{equation} The local submartingale $\left(Y_{t}\right)_{t%
\geq 0}$ then admits the following multiplicative decomposition:%
\begin{equation}  \label{multdecompformula}
Y_{t}=M_{t}C_{t}-1,
\end{equation}%
where $\left(M_{t}\right)_{t\geq 0}$ is a continuous local
martingale, which is strictly positive, with $M_{0}=1$ and where
$\left(C_{t}\right)_{t\geq 0}$ is an increasing continuous and
adapted process, with $C_{0}=1$. The decomposition is unique and the
processes $C$ and $M$ are given by the explicit formulae:
\begin{equation}  \label{partiecroissante}
C_{t}=\exp\left(\int_{0}^{t}\dfrac{d\ell_{s}}{1+Y_{s}}\right),
\end{equation}%
and
\begin{eqnarray}
M_{t} &=& \left(Y_{t}+1\right)\exp\left(-\int_{0}^{t}\dfrac{d\ell_{s}}{%
1+Y_{s}}\right)  \label{partiemartI} \\
&=& \exp\left(\int_{0}^{t}\dfrac{dm_{s}}{1+Y_{s}}-\dfrac{1}{2}\int_{0}^{t}%
\dfrac{d\langle m\rangle_{s}}{\left(1+Y_{s}\right)^{2}}\right).
\label{partiemartII}
\end{eqnarray}
\end{prop}

\begin{proof}
Assume first that such a representation exists; then from It\^{o}'s
formula:
\begin{equation*}
dY_{t}=dm_{t}+d\ell_{t}=C_{t}dM_{t}+M_{t}dC_{t};
\end{equation*}%
from this equality, we deduce that:
\begin{eqnarray*}
dm_{t} &=& C_{t}dM_{t} \\
d\ell_{t} &=& M_{t}dC_{t}.
\end{eqnarray*}%
Multiplying the second equality by $C_{t}$ yields:
\begin{equation*}
dC_{t}=C_{t}\dfrac{d\ell_{t}}{1+Y_{t}},
\end{equation*}
which is equivalent to:
\begin{equation*}
C_{t}=\exp\left(\int_{0}^{t}\dfrac{d\ell_{s}}{1+Y_{s}}\right).
\end{equation*}%
This establishes (\ref{partiecroissante}), and (\ref{partiemartI})
is an immediate consequence.To prove (\ref{partiemartII}), we notice
that:
\begin{equation*}
dM_{t}=M_{t}\dfrac{dm_{t}}{1+Y_{t}},
\end{equation*}%
which is equivalent to:
\begin{equation*}
M_{t}=\exp\left(\int_{0}^{t}\dfrac{dm_{s}}{1+Y_{s}}-\dfrac{1}{2}\int_{0}^{t}%
\dfrac{d\langle m\rangle_{s}}{\left(1+Y_{s}\right)^{2}}\right).
\end{equation*}
This completes the uniqueness part of the proof. The existence part
is immediate.
\end{proof}

One is often interested by some questions of integrability about
stochastic processes. The next proposition  gives some necessary and
sufficient conditions on $M$ and $C$ for $Y$ to be of class $(D)$.
Before stating it, we prove a general lemma.
\begin{lem}\label{ptitlem}
Let $\left(M_{t}\right)_{t\geq 0}$ be a nonnegative uniformly
integrable martingale, with $M_{0}=1$, and let
$\left(C_{t}\right)_{t\geq 0}$ be an increasing right continuous and
adapted process such that $C_{0}=1$. Then, we have:
$$\mathbb{E}\left(M_{\infty}C_{\infty}\right)=1+\mathbb{E}\left(\int_{0}^{\infty}M_{s}dC_{s}\right).$$
\end{lem}
\begin{proof}
We have:
$$\mathbb{E}\left(M_{\infty}C_{\infty}\right)=1+\mathbb{E}\left(M_{\infty}\int_{0}^{\infty}dC_{s}\right)=1+\mathbb{E}\left(\int_{0}^{\infty}M_{\infty}dC_{s}\right).$$
Now, the result of the lemma follows by replacing the constant
stochastic process  $\left(s,\omega\right)\rightarrow
M_{\infty}\left(\omega\right)$ with its optional projection
$\left(s,\omega\right)\rightarrow M_{s}\left(\omega\right)$.
\end{proof}
\begin{prop}\label{appartennaceclasseD}
Let $\left(Y_{t}\right)_{t\geq 0}$ be a continuous nonnegative local
submartingale such that $Y_{0}=0$, and consider its multiplicative
decomposition given by (\ref{multdecompformula}). Then, the
following are equivalent:
\begin{enumerate}
\item $\left(Y_{t}\right)_{t\geq 0}$ is of class $(D)$;
\item $\left(M_{t}\right)_{t\geq 0}$ is a uniformly integrable
martingale and
$$\mathbb{E}\left[M_{\infty}C_{\infty}\right]<\infty.$$
\item $\left(M_{t}\right)_{t\geq 0}$ is a uniformly integrable
martingale and
$$\mathbb{E}\left[\int_{0}^{\infty}M_{s}dC_{s}\right]<\infty.$$
\end{enumerate}
\end{prop}
\begin{proof}
$(1)\Rightarrow (2)$. Since $Y_{t}+1=M_{t}C_{t}$ and $C_{t}\geq1$,
we have:
$$M_{t}=\dfrac{Y_{t}+1}{C_{t}}\leq Y_{t}+1,$$and as $Y$ is of class
$(D)$, so is $M$. Furthermore,
$$M_{\infty}C_{\infty}=1+Y_{\infty}$$is integrable and this
completes the proof of $(2)$.

$(2)\Rightarrow (1)$. We have:
$$M_{t}C_{t}=C_{t}\mathbb{E}\left[M_{\infty}|\mathcal{F}_{t}\right]\leq\mathbb{E}\left[C_{\infty}M_{\infty}|\mathcal{F}_{t}\right],$$and
consequently, if $(2)$ is satisfied, then $Y$ is of class $(D)$.

$(2)\Leftrightarrow (3)$. This last equivalence is a consequence of
Lemma \ref{ptitlem}.
\end{proof}

\subsection{A multiplicative characterization of the class $\mathbf{(\Sigma)}$ and applications}
In this subsection, we shall see that Proposition
\ref{decompositionmultipl} takes a special form which characterizes
the local submartingales of the class $(\Sigma)$. The
characterization that follows should be compared with the
multiplicative characterization of the Az\'{e}ma's supermartingales
associated with honest times given in \cite{ashyordoob}.

Let $X$ be a local submartingale of class $(\Sigma)$ given by
(\ref{dcodexdbmey}). Since
$\left(dA_{t}\right)$ is carried by the zeros of $X$: $\left\{t:\;X_{t}=0%
\right\}$, (\ref{partiecroissante}) simply writes:
\begin{equation*}
C_{t}=\exp\left(A_{t}\right).
\end{equation*}
In fact, we still have a better result: denoting,
\begin{equation*}
I_{t}=\inf_{s\leq t}M_{s},
\end{equation*}%
we can prove that:%
\begin{equation*}
C_{t}=\dfrac{1}{I_{t}}.
\end{equation*}

\begin{prop}\label{decomultsigma}
Let $\left(X_{t}=N_{t}+A_{t},\;t\geq 0\right)$ be a nonnegative,
continuous local submartingale with $X_{0}=0$. Then the following
are equivalent:

\begin{enumerate}
\item $\left(X_{t},\;t\geq 0\right)$ is of class $(\Sigma)$, i.e. $\left(dA_{t}\right)$ is carried by the set: $\left\{t:\;X_{t}=0%
\right\}$.

\item There exists a unique strictly positive, continuous local martingale $%
\left(M_{t}\right)$, with $M_{0}=1$, such that:
\begin{equation}
X_{t}=\dfrac{M_{t}}{I_{t}}-1,
\end{equation}%
where
\begin{equation*}
I_{t}=\inf_{s\leq t}M_{s}.
\end{equation*}
The local martingale $\left(M_{t}\right)$ is given by:
\begin{equation}
M_{t}=\left(1+X_{t}\right)\exp\left(-A_{t}\right).
\end{equation}
\end{enumerate}
\end{prop}

\begin{proof}
$(1)\Rightarrow (2)$ Since $C_{t}=\exp\left(A_{t}\right)$, $dC_{t}$
is also carried by the set $\left\{t:\;X_{t}=0\right\}$. We also
have:
\begin{equation*}
C_{t}-\dfrac{1}{M_{t}}=\dfrac{X_{t}}{M_{t}}.
\end{equation*}
An application of Skorokhod's reflection lemma yields:
\begin{equation*}
C_{t}=\sup_{s\leq t}\dfrac{1}{M_{s}};
\end{equation*}%
that is:
\begin{equation*}
C_{t}=\dfrac{1}{I_{t}}.
\end{equation*}
$(2)\Rightarrow (1)$ Let $$X_{t}=\dfrac{M_{t}}{I_{t}}-1,$$with $M$ a
continuous and strictly positive local martingale. It is clear that
$X\geq0$. An application of It\^{o}'s formula further yields:
$$X_{t}=\int_{0}^{t}\dfrac{dM_{s}}{I_{s}}+\log\left(\dfrac{1}{I_{t}}\right).$$Consequently,
$X$ is of class $(\Sigma)$.
\end{proof}
\begin{rem}
In the proof of Proposition \ref{decomultsigma}, we have proved the
following fact:
$$X_{t}=\dfrac{M_{t}}{I_{t}}-1$$is a local submartingale of the
class $(\Sigma)$ and its Doob-Meyer decomposition is given by:
\begin{equation}\label{dbdesigma}
    X_{t}=\int_{0}^{t}\dfrac{dM_{s}}{I_{s}}+\log\left(\dfrac{1}{I_{t}}\right).
\end{equation}
\end{rem}
The above proposition underlines the fact that the local
submartingales of class $(\Sigma)$ are the smallest local
submartingales with a given martingale part $M$ in their
multiplicative decomposition. More precisely, let
$\left(M_{t}\right)$ be a strictly positive and continuous local
martingale with $M_{0}=1$ and denote $\mathcal{E}\left(M\right)$ the
set of all nonnegative local submartingales with the same martingale
part $M$ in the multiplicative decomposition
(\ref{multdecompformula}). Then the following holds:
\begin{cor}\label{ferquenceofvainsihing}
Let
$$Y^{\star}=\dfrac{M_{t}}{I_{t}}-1.$$Then,
$\left(Y^{\star}_{t}\right)$ is in $\mathcal{E}\left(M\right)$ and
it is the smallest element of $\mathcal{E}\left(M\right)$ in the
sense that:
$$\forall Y\in\mathcal{E}\left(M\right),\;Y^{\star}\leq
Y.$$Consequently, $\left(Y^{\star}_{t}\right)$ has more zeros than
any other local submartingale of $\mathcal{E}\left(M\right)$.
\end{cor}
\begin{proof}
It suffices to note that any element $Y\in\mathcal{E}\left(M\right)$
decomposes as $Y_{t}=M_{t}C_{t}-1$. Since $Y$ must be nonnegative,
we must have:
$$M_{t}\geq\dfrac{1}{C_{t}}.$$ But $\frac{1}{C_{t}}$ is decreasing,
hence we have:
$$\dfrac{1}{C_{t}}\leq I_{t},$$and this proves the Corollary.
\end{proof}
We can restate the above corollary as:
\begin{cor}
Let $Y$ be a nonnegative local submartingale, with $Y_{0}=0$. Then,
there exists a local submartingale $Y^{\star}$ of class $(\Sigma)$
such that:
$$Y^{\star}_{t}\leq Y_{t},\quad \mathrm{for\; all\;}t.$$
\end{cor}

This last property should be compared with a result of Az\'{e}ma and
Yor (\cite{azemayorzero}), according to which the set of zeros of a
process of the class $(\Sigma D)$ is a saturated set; we shall come
back on this point more precisely in the next subsection. For now,
let us give some examples and some other applications.
\begin{ex}
Let $\left(K_{t}\right)$ be a continuous local martingale, with
$K_{0}=0$. Then $\left(|K_{t}|\right)$ is of class $(\Sigma)$. From
Proposition \ref{decomultsigma}, and Tanaka's formula we have:
$$|K_{t}|=\exp\left(L_{t}\right)\exp\left(\int_{0}^{t}\dfrac{sgn\left(K_{s}\right)dK_{s}}{1+K_{s}}-\dfrac{1}{2}\int_{0}^{t}%
\dfrac{d\langle
K\rangle_{s}}{\left(1+K_{s}\right)^{2}}\right)-1,$$where
$\left(L_{t}\right)$ is the local time at $0$ of $K$. Similarly, we
have:
$$K_{t}^{+}=\exp\left(\dfrac{1}{2}L_{t}\right)\exp\left(\int_{0}^{t}\mathbf{1}_{\left(K_{s}>0\right)}\dfrac{dK_{s}}{1+K_{s}}-\dfrac{1}{2}\int_{0}^{t}%
\mathbf{1}_{\left(K_{s}>0\right)}\dfrac{d\langle
K\rangle_{s}}{\left(1+K_{s}\right)^{2}}\right)-1,$$
\end{ex}
\begin{ex}
With the same notations as above, if we define further:
$$\overline{K}_{t}=\sup_{s\leq t}K_{s},$$we have:
$$\overline{K}_{t}-K_{t}=\exp\left(\overline{K}_{t}\right)\exp\left(-\int_{0}^{t}\dfrac{dK_{s}}{1+K_{s}}-\dfrac{1}{2}\int_{0}^{t}%
\dfrac{d\langle K\rangle_{s}}{\left(1+K_{s}\right)^{2}}\right)-1.$$
\end{ex}
Now, let us give other applications of Proposition
\ref{decomultsigma}. The first one states that the local time at $0$
of a continuous local martingale is always equal to a function of
the infimum of some nonnegative local martingale. More precisely:
\begin{cor}
Let $\left(K_{t}\right)$ be a continuous local martingale, with
$K_{0}=0$, and let $\left(L_{t}\right)$ be its local time at $0$.
Then there exists a unique continuous and strictly positive local
martingale $\left(M_{t}\right)$, with $M_{0}=1$, such that:
$$L_{t}=\log\left(\dfrac{1}{I_{t}}\right),$$where $I_{t}=\inf_{s\leq
t}M_{s}$.
\end{cor}
\begin{proof}
From Tanaka's formula, we have
$$|K_{t}|=\int_{0}^{t}sgn\left(K_{s}\right)dK_{s}+L_{t}.$$But since
$|K|$ is of class $(\Sigma)$, from Proposition \ref{decomultsigma}
and formula (\ref{dbdesigma}), there exists a unique continuous and
strictly positive local martingale $\left(M_{t}\right)$, with
$M_{0}=1$, such that:
$$|K_{t}|=\int_{0}^{t}\dfrac{dM_{s}}{I_{s}}+\log\left(\dfrac{1}{I_{t}}\right).$$Now,
from the unicity of the additive decomposition, we obtain:
$$L_{t}=\log\left(\dfrac{1}{I_{t}}\right).$$
\end{proof}
\begin{cor}
Let $\left(K_{t}\right)$ be a continuous local martingale, and
define:
$$\overline{K}_{t}=\sup_{s\leq t}K_{s}.$$Then, there exist  two nonnegative and continuous local martingales
$\left(M_{t}^{(1)}\right)$ and $\left(M_{t}^{(2)}\right)$, with
$M_{0}^{(1)}=M_{0}^{(2)}=1$, such that:
\begin{eqnarray*}
  \left\{t:\;K_{t}=0\right\} &=& \left\{t:\;M_{t}^{(1)}=I_{t}^{(1)}\right\}; \\
  \left\{t:\;K_{t}=\overline{K}_{t}\right\} &=&
  \left\{t:\;M_{t}^{(2)}=I_{t}^{(2)}\right\}.
\end{eqnarray*}
\end{cor}
\begin{proof}
It suffices to note that
$\left\{t:\;K_{t}=0\right\}=\left\{t:\;|K_{t}|=0\right\}$ and then
apply Proposition \ref{decomultsigma} and formula (\ref{dbdesigma}).
\end{proof}
\begin{cor}\label{surmartazemasousmart}
Let $L$ be an honest time. Assume that $L$ avoids stopping times,
i.e. for all $\left(\mathcal{F}_{t}\right)$ stopping times $T$,
$\mathbb{P}\left(L=T\right)=0$. If furthermore
$\mathbb{P}\left(L\leq t|\mathcal{F}_{t}\right)$ is continuous (this
happens in particular when all $\left(\mathcal{F}_{t}\right)$
martingales are continuous), then there exists a unique continuous
and strictly positive local martingale $\left(M_{t}\right)$, with
$M_{0}=1$, such that:
$$\mathbb{P}\left(L\leq t|\mathcal{F}_{t}\right)=\dfrac{M_{t}}{I_{t}}-1=\int_{0}^{t}\dfrac{dM_{s}}{I_{s}}+\log\left(\dfrac{1}{I_{t}}\right).$$
Furthermore, we have:
$$M_{\infty}=2I_{\infty}.$$
\end{cor}
\begin{proof}
The result follows from the fact that $\left(\mathbb{P}\left(L\leq
t|\mathcal{F}_{t}\right)\right)$ is of class $(\Sigma)$ (see
\cite{azema}) and $\lim_{t\rightarrow\infty}\mathbb{P}\left(L\leq
t|\mathcal{F}_{t}\right)=1$.
\end{proof}
\begin{rem}
For the reader who is acquainted with the results of the general
theory of stochastic processes and in particular with the theory of
progressive enlargements of filtrations, it is easy to deduce from
the above corollary that in fact $L$ may be represented as:
$$L=\sup\left\{t:\;M_{t}=I_{t}\right\}=\inf\left\{t:M_{t}=I_{\infty}\right\}.$$Indeed,
from the results of Az\'{e}ma (\cite{azema}), $L$ is the last zero
of $\mathbb{P}\left(L\leq t|\mathcal{F}_{t}\right)$.
\end{rem}\bigskip

Now, let us give some necessary and sufficient conditions on $M$ for
$X$ to be of class $(\Sigma D)$. Let $\left(X_{t}\right)$ be of
class $(\Sigma)$: $$X_{t}=N_{t}+A_{t}=\dfrac{M_{t}}{I_{t}}-1.$$ Let
us also define
$$U_{t}=\int_{0}^{t}\dfrac{dM_{s}}{M_{s}}=\int_{0}^{t}\dfrac{dN_{s}}{1+X_{s}}.$$
\begin{prop}
Let $\left(X_{t}\right)$ be of class $(\Sigma)$. Then the following
properties are equivalent:
\begin{enumerate}
\item $\left(X_{t}\right)$ is of class $(\Sigma D)$;
\item $\left(M_{t}\right)$ is a uniformly integrable martingale and
$$\mathbb{E}\left(\dfrac{M_{\infty}}{I_{\infty}}\right)<\infty;$$
\item $\left(M_{t}\right)$ is a uniformly integrable martingale and
$$\mathbb{E}\left(\log\left(\dfrac{1}{I_{\infty}}\right)\right)<\infty;$$
\item $\left(M_{t}\right)$ is a uniformly integrable martingale and
$\left(U_{t}\right)$ is an $L^{2}$ bounded martingale.
\end{enumerate}
\end{prop}
\begin{proof}
The equivalence between $(1)$, $(2)$ and $(3)$ follows from
Proposition \ref{appartennaceclasseD} and the fact that here:
$$M_{s}dC_{s}=M_{s}\dfrac{\left(-dI_{s}\right)}{I_{s}^{2}}=-\dfrac{dI_{s}}{I_{s}}=d\left(\log\dfrac{1}{I_{s}}\right).$$
It now suffices to show for example $(3)\Leftrightarrow (4)$. First,
from It\^{o}'s formula, we note that:
\begin{equation}\label{foruleitolog}
\log\dfrac{1}{M_{t}}=-\int_{0}^{t}\dfrac{dM_{s}}{M_{s}}+\dfrac{1}{2}\int_{0}^{t}\dfrac{d\langle
M \rangle_{s}}{M_{s}^{2}}.
\end{equation}
Let us first show that $(4)\Rightarrow (3)$. From
(\ref{foruleitolog}), we have:
\begin{eqnarray*}
  \mathbb{E}\left[\sup_{s\leq t}|\log\dfrac{1}{M_{s}}|\right] &\leq& \mathbb{E}\left[\sup_{s\leq t}|\int_{0}^{s}\dfrac{dM_{u}}{M_{u}}|+\dfrac{1}{2}\int_{0}^{t}\dfrac{d\langle
M \rangle_{s}}{M_{s}^{2}}\right] \\
   &\leq& \mathbb{E}\left[\sup_{s\leq t}|U_{s}|\right]+\dfrac{1}{2}\mathbb{E}\left[\langle U\rangle_{t}\right] \\
   &\leq& C\mathbb{E}\left[\langle U\rangle_{t}^{1/2}\right]+\dfrac{1}{2}\mathbb{E}\left[\langle U\rangle_{t}\right] \\
   &\leq& C\left(\mathbb{E}\left[\langle
   U\rangle_{t}\right]\right)^{1/2}+\dfrac{1}{2}\mathbb{E}\left[\langle
   U\rangle_{t}\right],
\end{eqnarray*}where the third inequality is obtained by an
application of the Burkholder-Davis-Gundy inequalities and the last
inequality is a consequence of Jensen's inequality ($C$ stands for a
constant). This shows $(3)$.

Now, let us prove that $(3)\Rightarrow (4)$. Again, from
(\ref{foruleitolog}),and the Burkholder-Davis-Gundy inequalities, we
have, for any stopping time $T$:
$$\mathbb{E}\left[\langle
   U\rangle_{T}\right]\leq C\left(\left(\mathbb{E}\left[\langle
   U\rangle_{T}\right]\right)^{1/2}+\mathbb{E}\left[\log\dfrac{1}{I_{\infty}}\right]\right).$$Now,
   we can choose a sequence of stopping times $\left(T_{n}\right)$,
   such that $\lim_{n\rightarrow\infty}T_{n}=+\infty$ and $\mathbb{E}\left[\langle
   U\rangle_{T_{n}}\right]<\infty$. Dividing both sides of the above
   equality by $\left(\mathbb{E}\left[\langle
   U\rangle_{T_{n}}\right]\right)^{1/2}$ and letting
   $n\rightarrow\infty$, we obtain by an application of the monotone
   convergence theorem that:
   $$\left(\mathbb{E}\left[\langle
   U\rangle_{\infty}\right]\right)^{1/2}\leq C\left(1+\dfrac{1}{\left(\mathbb{E}\left[\langle
   U\rangle_{\infty}\right]\right)^{1/2}}\mathbb{E}\left[\log\dfrac{1}{I_{\infty}}\right]\right),$$and
   this shows that
   $$\mathbb{E}\left[\langle
   U\rangle_{\infty}\right]<\infty.$$
\end{proof}
\begin{rem}
In the course of the proof of the equivalence between assertions
$(3)$ and $(4)$, we did not use the fact that $\left(M_{t}\right)$
is a uniformly integrable martingale. In fact, we have proved the
following fact:

If $M_{t}\equiv \exp\left(U_{t}-\dfrac{1}{2}\langle
U\rangle_{t}\right)$, where $\left(U_{t}\right)$ is a continuous
local martingale, such that $U_{0}=0$, then we have:
$$\mathbb{E}\left[\log\dfrac{1}{I_{\infty}}\right]<\infty\Leftrightarrow \mathbb{E}\left[\langle
   U\rangle_{\infty}\right]<\infty.$$One moment of attention shows
   that $$\log\dfrac{1}{I_{\infty}}=-\inf_{t\geq0}\left(U_{t}-\dfrac{1}{2}\langle
U\rangle_{t}\right),$$and hence:
$$\mathbb{E}\left[\langle
   U\rangle_{\infty}\right]<\infty\Leftrightarrow \mathbb{E}\left[\inf_{t\geq0}\left(U_{t}-\dfrac{1}{2}\langle
U\rangle_{t}\right)\right]>-\infty.$$
\end{rem}
\begin{rem}
The property of being uniformly integrable for $M_{t}=
\exp\left(U_{t}-\dfrac{1}{2}\langle U\rangle_{t}\right)$ is not
necessary, nor sufficient to have $\mathbb{E}\left[\langle
   U\rangle_{\infty}\right]<\infty$. Indeed, from Dubins-Schwarz
   theorem,
   $$U_{t}=B_{\langle U\rangle_{t}},$$where $B$ is a standard
   Brownian Motion. Hence, noting $T=\langle U\rangle_{\infty}$, it
   is equivalent to show that $$M_{t}\equiv\exp\left(B_{t\wedge T}-\dfrac{1}{2}t\wedge
   T\right)$$ may be uniformly integrable without having
   $\mathbb{E}\left(T\right)<\infty$, and vice versa.

   Let $T\equiv T_{a}=\inf\left\{t:\;B_{t}\geq a\right\}$, with
   $a>0$; then the local martingale $\left(M_{t}\right)$ defined
   above is bounded and hence uniformly integrable; however, it is
   well known that $\mathbb{E}\left(T_{a}\right)=\infty$.

   Now, conversely, let $$\sigma_{b}\equiv\inf\left\{t:\;B_{t}+bt=
   1\right\},$$ with $b>0$. From \cite{revuzyor} (exercise 1.31,
   p.335),
   $$\mathbb{E}\left(\exp\left(B_{\sigma_{b}}-\dfrac{1}{2}\sigma_{b}\right)\right)<1,$$and
   consequently
   $M_{t}=\exp\left(B_{t\wedge\sigma_{b}}-\dfrac{1}{2}t\wedge\sigma_{b}\right)$
   is not uniformly integrable. Nevertheless, $\sigma_{b}$ has
   exponential moments (\cite{revuzyor}, exercise 1.21, p.334).
\end{rem}
\subsection{Saturated sets}
The aim of this paragraph is to understand better Corollary
\ref{ferquenceofvainsihing} with the help of some results of
Az\'{e}ma, Meyer and Yor \cite{azemameyeryor,azemayorzero}.
Corollary \ref{ferquenceofvainsihing} tells us that the local
submartingales of the class $(\Sigma)$ vanish more than any other
local submartingale with the same local martingale component in its
multiplicative decomposition. Az\'{e}ma and Yor, in
\cite{azemayorzero}, have obtained a similar characterization, which
in some sense complements Corollary \ref{ferquenceofvainsihing}:
they proved that the set of zeros of submartingales of class
$(\Sigma D)$, whose terminal value is almost surely nonzero, is
saturated, that is every predictable set which is at the left of
such a random set is in fact contained into it. The notion of
saturated set was first introduced by P.A. Meyer (see
\cite{delmaismey,azemameyeryor}) and it is probably his last
contribution to the general theory of stochastic processes. It was
then used by Az\'{e}ma and Yor \cite{azemayorzero} in the study of
the set of zeros of continuous martingales and by Yor in an attempt
to prove Barlow's conjecture on Brownian filtrations (\cite{zurich},
and for a proof of Barlow's conjecture, see \cite{beksy}). First, we
give the definition of a saturated set, and then we cite a theorem
of Az\'{e}ma and Yor. Eventually, we give our own representation of
saturated sets, and then we conclude giving several nontrivial
examples of such sets.

In the sequel, for simplicity, we always assume that the filtration
$\left(\mathcal{F}_{t}\right)$ is such that every
$\left(\mathcal{F}_{t}\right)$ martingale is continuous. In
particular, the predictable and optional sigma fields are equal. We
assume we are given a predictable set $H$ whose end
$$g=\sup\left\{t:\;\left(t,\omega\right)\in H\right\}$$ is almost
surely finite. We shall also say that a predictable random set $E$
is at the left of $g$ (or $H$) if the end $L$ of $E$ satisfies:
$L\leq g$. There are several equivalent definitions of saturated
sets in the literature which we have collected below (see
\cite{delmaismey} p.135, \cite{azemameyeryor,azemayorzero},
\cite{zurich} p.108):
\begin{defn}\label{defsatur}
The following definitions of saturated sets are equivalent.
\begin{enumerate}
\item (\cite{zurich} p.108) The predictable set $H$ is said to be saturated if for every end $L$ of a predictable set $L\leq
g$ one has:
$$\left(L\left(\omega\right),\omega\right)\in H.$$
\item (\cite{azemayorzero}) The predictable set $H$ is saturated if
every predictable set which is at the left of $H$ is in fact
contained in $H$.
\item (\cite{delmaismey} p.135) The predictable set $H$ is saturated
if $$H=\left\{t:\;\mathbb{P}\left(g\leq
t|\mathcal{F}_{t}\right)=0\right\},$$or in the very imaged
terminology of Strasbourg, $H$ is equal to its predictable shadow.
\end{enumerate}
\end{defn} Assume
that $g$ avoids $\left(\mathcal{F}_{t}\right)$ stopping times, which
is equivalent (with our assumptions on
$\left(\mathcal{F}_{t}\right)$ ) to the fact that
$\left(\mathbb{P}\left(g\leq t|\mathcal{F}_{t}\right)\right)$ is
continuous. Consequently, from well known results by Az\'{e}ma,
$\left(\mathbb{P}\left(g\leq t|\mathcal{F}_{t}\right)\right)$ is of
class $(\Sigma D)$. Now we can state the result of Az\'{e}ma and Yor
which complements our characterization of local submartingales of
class $(\Sigma)$.
\begin{prop}[Az\'{e}ma-Yor \cite{azemayorzero}]
Let $\left(X_{t}\right)$ be submartingale of class $(\Sigma D)$,
such that $$\mathbb{P}\left(X_{\infty}=0\right)=0.$$Define
$$H=\left\{t:\;X_{t}=0\right\}.$$ Then $H$ is a saturated set.
\end{prop}
\begin{rem}\label{remapredef}
The characterization of Corollary \ref{ferquenceofvainsihing} on the
frequency of vanishing of the local submartingales of class
$(\Sigma)$ did not require integrability and limit conditions. In
the above proposition, the integrability condition is essential.
Indeed, take $$X_{t}=|B_{t\wedge T_{1}}|,$$where
$T_{1}=\inf\left\{t:B_{t}=1\right\}$. Here $$H=\left\{t\leq
T_{1}:\;B_{t}=0\right\}.$$ Now let us consider
$$L=\sup\left\{t\leq T_{1}:\;B_{t}=I_{t}\right\},$$where as usual
$I_{t}=\inf_{u\leq t}B_{u}$. Then $$X_{L}=|B_{L}|>0,$$and $$L\leq
g,$$ where $g=\sup\left\{t\leq T_{1}:\;B_{t}=0\right\}$.
Consequently, $H$ is not saturated.
\end{rem}
Now, let us give some other characterizations for saturated sets.
\begin{prop}\label{propvingt}
Let $H$ be a saturated set whose end $g$ avoids
$\left(\mathcal{F}_{t}\right)$ stopping times. Then, there exist two
nonnegative and continuous local martingale $\left(M_{t}\right)$ and
$\left(N_{t}\right)$, with $M_{0}=1,N_{0}=1$ and
$M_{\infty}=2I_{\infty}$ (where $I_{t}=\inf_{s\leq t}M_{s}$),
$N_{\infty}=0$, such that $H$ may be represented as:
$$H=\left\{t:\;M_{t}=I_{t}\right\}=\left\{t:\;N_{t}=\overline{N}_{t}\right\},$$where
$\overline{N}_{t}=\sup_{u\leq t}N_{u}$. Moreover, the local
martingales $\left(M_{t}\right)$ and $\left(N_{t}\right)$ are
unique.
\end{prop}
\begin{proof}
This is a consequence of Definition \ref{defsatur} (3) and Corollary
\ref{surmartazemasousmart} for the part with $M$. The part with $N$
follows from the multiplicative representation of
$\left(\mathbb{P}\left(g\leq t|\mathcal{F}_{t}\right)\right)$ proved
in \cite{ashyordoob}: there exists a (unique) continuous local
martingale $N$, with $N_{0}=1$ and $N_{\infty}=0$, such that
$$\mathbb{P}\left(g\leq t|\mathcal{F}_{t}\right)=1-\dfrac{N_{t}}{\overline{N}_{t}}.$$
\end{proof}

We end the paper with some non trivial examples of saturated sets
one may encounter in martingale theory or in the theory of diffusion
processes.
\begin{ex}
Let $\left(M_{t}\right)$ be a continuous and nonnegative local
martingale, such that $M_{0}=x>0$ and
$\lim_{t\rightarrow\infty}M_{t}=0$. Define, for $0<a\leq x$:
$$g_{a}=\sup\left\{t\geq0:\;M_{t}=a\right\}.$$ It is proved in
\cite{AshkanremII} that:
$$\mathbb{P}\left(g\leq
t|\mathcal{F}_{t}\right)=1-\left(\dfrac{M_{t}}{a}\right)\wedge1.$$Consequently,
the set
$$H=\left\{t:\;M_{t}\geq a\right\}$$is saturated.

We already noted in Remark \ref{remapredef} that the set
$\left\{t\leq T_{1}:\;B_{t}=0\right\}$ is not saturated. Now, if we
apply the previous result to the martingale $M_{t}=1-B_{t\wedge
T_{1}}$, we obtain that the set:
$$H=\left\{t\leq T_{1}:\;B_{t}\leq0\right\}$$ is saturated.

We can also apply the above result to $\left(R_{t}\right)$, a
transient diffusion with values in $\left[0,\infty\right)$, which
has $\left\{0\right\}$ as entrance boundary. Let $s$ be a scale
function for $R$, which we can choose such that:
$$s\left(0\right)=-\infty, \text{ and } s\left(\infty\right)=0.$$
Then, under the law $\mathbb{P}_{x}$, for any $x>0$, the local
martingale $\left(M_{t}=-s\left(R_{t}\right)\right)$ satisfies the
conditions above and for $0\leq x \leq y$, we
have:$$\mathbb{P}_{x}\left(g_{y}>
t|\mathcal{F}_{t}\right)=\dfrac{s\left(R_{t}\right)}{s\left(y\right)}\wedge
1,$$where$$g_{y}=\sup\left\{t:\; R_{t}=y\right\}.$$By Definition
2.17 (3), the set
$$H=\left\{t:\;R_{t}\leq y\right\}$$is saturated.
\end{ex}
\begin{ex}
Let $N_{t}=1-B_{t\wedge T_{1}}$. Then
$\overline{N}_{t}=1-\inf_{u\leq t\wedge T_{1}}B_{u}=1-I_{t\wedge
T_{1}}$ and $N_{\infty}=0$. Consequently, from Proposition
\ref{propvingt}, the random set
$$H=\left\{t\leq T_{1}:\;B_{t}=I_{t}\right\}$$is saturated.

Now, we consider again a transient diffusion $\left(R_{t}\right)$,
with values in $\left[0,\infty\right)$, which has $\left\{0\right\}$
as entrance boundary. The scale function $s$ is again chosen such
that:
$$s\left(0\right)=-\infty, \text{ and } s\left(\infty\right)=0.$$
Then, under the law $\mathbf{P}_{x}$, for any $x>0$, the local
martingale
$\left(N_{t}=\dfrac{s\left(R_{t}\right)}{s\left(x\right)},\;t\geq
0\right)$ satisfies the conditions of Proposition \ref{propvingt},
and we have:$$\mathbf{P}_{x}\left(g>
t|\mathcal{F}_{t}\right)=\dfrac{s\left(R_{t}\right)}{s\left(I_{t}\right)}\,$$where
$$g=\sup\left\{t:\; R_{t}=I_{t}\right\},$$and
$$I_{t}=\inf_{s\leq t}R_{s}.$$Hence, the set $$H=\left\{t\geq0:\;R_{t}=I_{t}\right\}$$is saturated.
\end{ex}
\section*{Acknowledgements}
I am deeply indebted to Marc Yor whose ideas and help were essential
for the development of this paper: he shared very kindly with me
some unpublished notes about saturated sets.

\end{document}